\documentclass[a4paper,11pt]{article}
\usepackage[bookmarks=true,
              bookmarksnumbered=true, breaklinks=true,
              pdfstartview=FitH, hyperfigures=false,
              plainpages=false, naturalnames=true,
              colorlinks=true,%pagebackref=true,
              pdfpagelabels]{hyperref}
\usepackage{amssymb,amscd,amsthm,amsmath,color,fancyhdr,mathrsfs,stackrel}
\usepackage{enumerate,enumitem}
\usepackage[all]{xy}
\usepackage{authblk}
\usepackage[left=1.5cm, right=1.5cm, top=2.5cm,
bottom=2.5cm]{geometry}
\usepackage{mathtools}
\usepackage{float}
\usepackage{lineno}
\usepackage{color,soul}
\usepackage{indentfirst}

\newtheorem{thm}{Theorem}[section]
\newtheorem{prop}[thm]{Proposition}
\newtheorem{lem}[thm]{Lemma}

\newtheorem{cor}[thm]{Corollary}

\theoremstyle{definition}
\newtheorem{defn}[thm]{Definition}
\newtheorem*{prf}{Proof of Proposition \ref{prop:51}}

\def\wtilde{\widetilde}
\def\inc{\mathbin{\widetilde\in}}
\def\incs{\mathbin{\widetilde\subset}}

\def\cIm{\operatorname{cIm}}

\def\cGr{\mathsf{cGr}}

\def\Sets{\mathsf{Sets}}

\def\cKer{\textbf{cKer}}

\def\cCrsM{\mathsf{cssc\text{-}XMod}}
\def\Catgr{\mathsf{{Catgr}}}
\def\qed{\hfill $\Box$}
\def\Ker{\operatorname{Ker}}
\def\cKer{\operatorname{cKer}}
\def\NSCM/(A,B,\mu){\mathsf{NSCM/(A,B,\mu)}}
\def\NSGGd/G{\mathsf{NSGGd/G}}

\def\epsilon{\varepsilon}
\def\wtilde{\widetilde}

\def\C{\mathsf{C}}
\def\Str{\textbf{Star}_{\mathsf{C}} 0}
\def\Strp{\textbf{Star}_{\mathsf{C'}} 0}
\def\Strb{\textbf{Star}_{\mathsf{\overline{C}}} 0}

\numberwithin{equation}{section}

\title{\bf Equivalence of the categories of categorical groups and cssc-crossed modules}

\author[1]{Tamar Datuashvili\thanks{\textbf{Corresponding author:} Tamar Datuashvili (e-mail: tamar.datu@gmail.com)}}
\affil[1]{\small{A. Razmadze Mathematical Institute of I. Javakhishvili Tbilisi State University, 2, M. Aleksidze II Lane, 0193, Tbilisi, GEORGIA}}
\author[2]{Osman Mucuk\thanks{O. Mucuk (e-mail : mucuk@erciyes.edu.tr)}}
\affil[2]{\small{Department of Mathematics, Erciyes University, 38039, Kayseri, TURKEY}}
\author[2]{Nazmiye Alemdar\thanks{N. Alemdar (e-mail : nakari@erciyes.edu.tr)}}
\author[3]{Tun\c{c}ar \c{S}ahan\thanks{T. \c{S}ahan (e-mail : tuncarsahan@gmail.com)}}
\affil[3]{\small{Department of Mathematics, Aksaray University, 68100, Aksaray, TURKEY}}

\date{\textit{Dedicated to the memory of Professor Ronald Brown}}

\begin{document}
%\linenumbers	
\maketitle

\begin{abstract}
In \cite{Dat2020} we proved that any categorical group defines a c-crossed module, which is a cssc-crossed module defined in the same paper. In \cite{Dat2024} we constructed a categorical group for any cssc-crossed module. In the presented paper we prove that these correspondences define functors between the corresponding categories, which realize an equivalence of the categories of categorical groups and cssc-crossed modules.
\end{abstract}	

\section{Introduction}

This paper is a continuation of our previous works \cite{Dat2020,Dat2024}. In the first one we introduced notions of c-groups, which is a group up to congruence relation, of c-crossed module and cssc-crossed module, which is a connected, strict and special c-crossed module, defined in the same work.

By a categorical group we mean a coherent 2-group in the sense of Baez and Lauda \cite{baez-lauda-2-groups}. In the same work it is proved that any categorical group gives rise to a cssc-crossed module. In \cite{Dat2024} for any csssc-crossed module we constructed a category equipped with a structure and proved that this is a categorical group. It is a well known fact that a categorical group is equivalent to a strict categorical group \cite{baez-lauda-2-groups,JoyalStreet,Sinh 2}, but we do not have equivalence between the corresponding categories.

Our goal was to obtain for categorical groups an analogous description in terms of certain crossed module type objects as we have for $\mathcal{G}-$groupoids obtained by Brown and Spencer \cite{BS1}, which are strict categorical groups, or equivalently, group-groupoids or internal categories in the category of groups. In this paper we prove that the above noted constructions define the functors between the categories of categorical groups and cssc-crossed modules and they realize an equivalence of these categories.

We hope that this result will give a chance to consider for categorical groups problems analogous to these considered and solved in the case of strict categorical groups in terms of group-groupoids and internal categories in \cite{BM1,Dat1,Dat2,Dat3,Dat4,Mucuk2015}.

\section{Preliminary results and definitions}

Since the paper is a continuation of the previous papers, in this section as a preparation we summarise some preliminaries from \cite{Dat2020, Dat2024}.

%Since the paper is a continuation of the previous paper, in this  section we remind some necessarily preliminaries from \cite{Dat2020}.

\subsection{Categorical groups}

The following definition of a monoidal category goes back to Mac Lane \cite{Mac71}.

A \textbf{monoidal category} is a category $\mathsf{C}=(C_0, C_1, d_0, d_1, i, m)$ equipped with a bifunctor $+\colon\mathsf{C}\times\mathsf{C}\rightarrow\mathsf{C}$ called the monoidal sum, an object $0$ called the zero object, and the family of three natural isomorphisms  $\alpha=\alpha_{x,y,z}\colon (x+y)+z\stackrel{\approx}{\rightarrow} x+(y+z)$, $\lambda_x\colon 0+x\stackrel{\approx}{\rightarrow} x$ and $\rho_x\colon x+0\stackrel{\approx}{\rightarrow} x$, for all $x,y,z\in C_0$, such that the following diagrams commute for all $x,y,z,t\in C_0$:
\[\xymatrix{
	((x+y)+z)+t  \ar[d]_{\alpha+1} \ar[r]^-{\alpha} & (x+y)+(z+t)  \ar[r]^-{\alpha} & x+(y+(z+t)) \\
	(x+(y+z))+t  \ar[rr]_-\alpha & & x+((y+z)+t) \ar[u]_{1+\alpha} }\]
\[\xymatrix{
	(x+0)+y  \ar[r]^-{\alpha} \ar[d]_-{\rho+1} &  x+(0+y)  \ar[d]^-{1+\lambda} \\
	x+y  \ar@{=}[r] & x+y }\]
$\lambda_0=\rho_0\colon 0+0\approx 0$.

Moreover, all diagrams involving $\alpha,\lambda$, and $\rho$ commute. A monoidal groupoid is a monoidal category in which every morphism is invertible.

In this definition we use the term monoidal sum and denote it as $+$, instead of monoidal product, used in the original definition, and write the operation additively.

It follows from the definition that $1_0+f\approx f+1_0 \approx f$, for any morphism $f$. A monoidal category is said to be strict if the natural isomorphisms $\alpha, \lambda$, and $\rho$ are identities.

For any two monoidal categories $\mathsf{C}=(C,+,0,\alpha,\lambda,\rho)$ and $\mathsf{C'}=(C',+',0',\alpha',\lambda',\rho')$, a functor $T:C\rightarrow C'$ satisfying $T(x+y)=Tx+'Ty$, $T(f+g)=Tf+'Tg, T0=0'$, $T\alpha_{x,y,z}=\alpha'_{Tx,Ty,Tz}$, $T\lambda_x=\lambda'_{Tx}$, $T\rho_x=\rho'_{Tx}$ for all objects $x,y,z$ and morphisms $f$ and $g$ is called a (strict) morphism of monoidal categories \[T:(C,+,0,\alpha,\lambda,\rho)\rightarrow (C',+',0',\alpha',\lambda',\rho').\]

Let $x$ be an object in a monoidal category $\mathsf{C}$. If there is an object $y\in C_0$ such that $x+y\approx 0$ and $y+x\approx 0$ then $y$ is called an inverse of $x$. It is a well known fact that if any object has a one-sided inverse in a monoidal category, then any object is invertible \cite{baez-lauda-2-groups,JoyalStreet}.

A monoidal groupoid $\mathsf{C}=(C_0, C_1, d_0, d_1, i, m)$ is called a \textit{categorical group} if for every object $x\in C_0$ there is an object $-x\in C_0$ with a family of natural isomorphisms	$\epsilon_x\colon -x+x\approx 0$ and $\delta_x\colon x+(-x)\approx 0$ such that the following diagrams are commutative:
\[\xymatrix@C=2.3pc{
	0 + x \ar[r]^<<<<<<{\delta^{-1}_x + 1_x} \ar[d]_{\lambda _x}
	& (x + (-x)) + x \ar[r]^<<<<<{\alpha_{x, -x ,x}}
	& x+( -x + x) \ar[d]^{1_x + \epsilon_x} \\
	x \ar[rr]_{\rho ^{-1}_x}
	& & x + 0  }
\]

\[
\xymatrix@C2.3pc{
	-x + 0
	\ar[r]^<<<<<<{1_x + \delta^{-1}_x}
	\ar[d]_{\rho _{-x}}
	& -x + (x + (-x))
	\ar[r]^{\alpha^{-1}_{-x, x, -x}}
	& (-x + x)+ (-x)
	\ar[d]^{\epsilon_x + 1_{-x}} \\
	-x
	\ar[rr]_{\lambda^{-1}_{-x}}
	&& 0 + (-x).   }
\]
%$(1_x + \epsilon_x)\alpha_{x, -x ,x}(\delta^{-1}_x + 1)=\rho ^{-1}_x\lambda _x$ and $(\epsilon_x + 1_{-x})\alpha^{-1}_{-x, x, -x}(1 + \delta^{-1}_x)=\lambda^{-1}_{-x}\rho_{-x}$

It is important and  a well-known fact that the definition of a categorical group implies that for any morphism $f:x\rightarrow x' \in C_1$ there is a morphism $-f\colon -x\rightarrow -x'$ with natural isomorphisms $-f+f\approx 0$ and $f+ (-f)\approx 0$, where the morphism $0$ is $1_0$ (see e.g. \cite{Sinh 1}).
The natural isomorphisms $\alpha,\lambda$, $\rho, \epsilon, \delta$, identity transformation $1_ \mathsf{C}\rightarrow 1_\mathsf{C}$, their compositions and sums are called \emph{special isomorphisms} in \cite{Dat2020}. Categorical group defined above is coherent \cite{baez-lauda-2-groups,Laplaza}, which means that all diagrams commute involving special isomorphisms. See \cite[Chapter VII Section 2]{Mac71} for coherence of monoidal categories.

The functorial properties of addition $+$ implies that in a categorical group we have $-1_x=1_{-x}$, for any $x\in C_0$.

The following Lemma is an analogue of the result obtained for internal categories by Porter \cite{Por}. It was not stated and proved in the previous papers.

\begin{lem}\label{lem:compadd}
Composition in any categorical group $\mathsf{C}=(C_0, C_1, d_0, d_1, i, m)$ is defined up to special isomorphism by the addition operation
\begin{equation*}
f\circ g\approx(f-i d_0(f))+g
\end{equation*}
for any composable arrows $f,g\in \mathsf{C}$.
\end{lem}
\begin{proof}
Since the addition operation ``+'' is a functor $\mathsf{C}\times \mathsf{C}\rightarrow \mathsf{C}$ like in the case of internal categories \cite{Por} we have the interchange law
\begin{equation*}
(f+h)\circ (g+t) = f\circ g + h\circ t,
\end{equation*}
whenever each side has a sense. From which easily follows the equality
\begin{equation*}
f\circ g - g = f - i d_0(f),
\end{equation*}
which implies the desired isomorphism.
\end{proof}

Since an isomorphism between morphisms $\theta :f\thickapprox g$ means that there exist isomorphisms $\theta_i:d_i(f)\rightarrow d_i(g), i=0,1$ with $\theta_1 f=g\theta_0$, the naturality property of special isomorphisms implies that there exist special isomorphisms between the morphisms in $C_1$. But if $\theta_i, i=0,1$ are special isomorphisms, it does not imply that $\theta$ is a special isomorphism; in this case $\theta$ is called a \emph{weak special isomorphism}. It is obvious that a special isomorphism between the morphisms in $C_1$ is a weak special isomorphism. Note that if $f\approx f'$ is a weak special isomorphism, then the coherence property implies that $f'$ is the unique morphism weakly specially isomorphic to $f$ with the same domain and codomain objects as $f'$.

In \cite{Dat2020}, we defined (strict) morphisms between categorical groups, which satisfy conditions of (strict) morphisms of monoidal categories. Note that this definition implies: $T(-x)=-{T(x)}$ and $T(-f)=-{T(f)}$, for any object $x$ and arrow $f$ in a categorical group. Categorical groups form a category with (strict) morphisms between them. 	For any categorical group $\mathsf{C}=(C_0, C_1, d_0, d_1, i, m)$ denote $\Ker d_0=\{f\in C_1\mid d_0(f)\approx 0\}$  and $\Ker d_1=\{f\in C_1\mid d_1(f)\approx 0\}$.

\begin{lem}\label{comm}
	Let $\mathsf{C}=(C_0, C_1, d_0, d_1, i, m)$ be a categorical group. For any $f\in \Ker d_1$ and $g\in \Ker d_0$ we have a weak special isomorphism $f+g\approx g+f$.
\end{lem}

\subsection{Groups up to congruence relation}

Now we recall the definition of a group up to congruence relation or briefly a c-group. Let the pair $X_R$ denotes a set $X$ with an equivalence relation $R\subseteq X\times X$. These kind of objects form a category, denoted by $\wtilde{\Sets}$, where the morphisms are functions $f\colon X_R\rightarrow Y_S$, such that $f(x)\sim_S f(y)$, whenever $x\sim_R y$.

Product in this category consists of the cartesian product of the sets and the usual product of the equivalence relations.

\begin{defn}\label{defcgr}
	A \textit{c-group} is an object $G_R$ in $\wtilde{\Sets}$ with a morphism $\mathfrak{m}\in \wtilde{\Sets}((G\times G)_{R\times R},G_R)$, denoted by $\mathfrak{m}(a,b)=a+b$, for any $a,b\in G$, satisfying the following conditions:
	%
	%		Let $G_R$ be an object in $\wtilde{\Sets}$ and
	%		\[ \begin{array}{cccl}
		%			m \colon &G\times G& \longrightarrow & G\\
		%			&(a,b)    & \longmapsto     & a+b
		%		\end{array}\]
	%		a morphism in $\wtilde{\Sets}$, i.e.,  $m\in \wtilde{\Sets}((G\times G)_{R\times R},G_R)$. $G_R$ is called a {\em c-group} if the following axioms are satisfied.
	\begin{enumerate}[label=(\roman{*}), leftmargin=1cm]
		\item\label{def:cgri} $a+(b+c)\sim_R(a+b)+c$, for all $a,b,c\in G$;
		\item\label{def:cgrii}  there exists an element $0\in G$ such that $a+0\sim_R a\sim_R 0+a$, for all $a\in G;$
		\item\label{def:cgriii} for each $a\in G$ there exists an element $-a$ such that  $a+(-a)\sim_R 0$ and $-a+a\sim_R0$.
	\end{enumerate}
\end{defn}

In a c-group $G_R$, the given element $0\in G$ is called  {\em zero element} and for any $a\in G$ the given element $-a\in G$ is called the {\em inverse} of $a$. The congruences involved in the conditions \ref{def:cgri}--\ref{def:cgriii} of the definition, their compositions and sums are called \emph{special congruences}.

Every group $G$ can be considered as a c-group with $R=\{(a,a)\mid a\in G\}$, i.e. $R$ is the equality $(=)$ relation. See \cite{Dat2020} for the properties and examples of c-groups.

For any two c-groups $G_R$ and  $H_S$, a morphism $f\colon G_R\rightarrow H_S$ in $\wtilde{\Sets}$ is called a \textit{c-group morphism} if $f(a+b)=f(a)+f(b)$,  for any $a,b\in G$. It is easy to see that a morphism between c-groups preserves congruences between elements; moreover $f(0)\sim 0$ and $f(-a)\sim -f(a)$, for any $a\in G$. Hence, a morphism between c-groups carries special congruences to special congruences between pairs of elements.

%A category can be constructed with objects as c-groups and with morphisms as c-group morphisms. This category is denoted by $\cGr$.

A category with objects as c-groups, with morphisms as c-group morphisms, with obvious composition of morphisms and identity morphisms of each object will be denoted by $\cGr$.

For a c-group morphism $f\colon G_R\rightarrow H_S$, the subset $\cKer f=\{a\in G_R \mid f(a)\sim_S 0_H \}$ of $G_R$ is called the \textit{c-kernel} and the subset $\cIm f=\{b\in H_S \mid \exists a\in G_R,  f(a)\sim_S b\}$ of $H_S$ is said to be the \textit{c-image} of $f$.

Let $H$ be a subset of the underlying set $G$ of a c-group $G_R$ and let $S=R\cap(H\times H)$. In this case, it is easy to see that $S$ is an equivalence relation on $H$. If $H_S$ is a c-group with the operation induced from $G_R$ then $H_S$ is called a c-subgroup of $G_R$. Note that $\cKer f$ is a c-subgroup of $G_R$. In particular, $\cKer d_0$, for a categorical group $\mathsf{C}=(C_0, C_1, d_0, d_1, i, m)$, is a c-subgroup of $C_1$ with the congruence relation on $\cKer d_0$ induced by the isomorphisms in $C_1$.

Let $G_R$ be a c-group and let $H$ be a subset of $G$. If for an element $a\in G$ there exists an element $b\in H$ such that $a\sim_R b$ then we write $a\inc H$. If $H$ and $H'$ are two subsets of $G_R$, then we write $H\incs H'$ if for any $h\in H$ we have $h \inc H'$. If $H\incs H'$ and $H'\incs H$, then we write $H\sim H'$.

\begin{defn}
	Let $H_S$ be a c-subgroup of a c-group $G_R$. Then $H_S$ is called
	\begin{enumerate}[label=(\roman{*}), leftmargin=1cm]
		\item \textit{normal} if $g+h-g\inc H_S$, for any $h\in H_S$ and $g\in G$;
		\item \textit{perfect} if $g\inc H$ implies $g\in H$, for any $g\in G$.
	\end{enumerate}
\end{defn}

One can see that, for a c-group morphism $f\colon G_R\rightarrow H_S$, $\cKer f$ is a perfect and normal c-subgroup of $G_R$, and $\cIm f$ is a perfect c-subgroup of $H_S$.

\begin{defn}
	A c-group $G_R$ is called \emph{connected} if $g\sim g'$ for any $ g, g'\in G$.
\end{defn}

%Quotient object $G/H$ of c-groups where $H$ is a normal c-subgroup of a c-group $G$ is described in \cite{Dat2020} as in the following way: Consider the classes  $\{g+H\mid g\in G\}$. If $(g+H) \cap (g'+H)\neq \emptyset$, then we obtain $-g+g'\inc H$, which implies that $g+H\sim g'+H$. Now consider the set of these classes $\{\mathrm{cl}(g+H)\mid g\in G\}$, where $\mathrm{cl}(g+H)=\cup\{x\in G\mid x\inc g+H\}$. Define $G/H=\{\mathrm{cl}(g+H)\mid g\in G\}$. An addition operation in this set is defined by $\mathrm{cl}(g+H)+\mathrm{cl}(g'+H)=\mathrm{cl}((g+g')+H)$, for any $g,g'\in G$. This operation is well defined, it is associative and we have the unit element $\mathrm{cl}(0+H)$. Actually the constructed object is a group, the congruence  relation on $G/H$ is the equality $``="$. We have an obvious surjective morphism $p: G \rightarrow G/H$.
%	
%	\begin{lem}\cite{Dat2020} (i) If $G$ is a c-group and $H$ is a normal c-subgroup in $G$, then for any group $G'$ and c-group morphism $f:G\rightarrow G'$, if $f(h)=0$ for any $h\in H$, there exists a unique morphism $\theta :G/H \rightarrow G'$, in $\cGr$  such that $\theta p=f$.\newline
	%		(ii)  If $H$ is a perfect normal c-subgroup in $G$, then $H=\cKer p$.
	%	\end{lem}

\subsection{Actions and crossed modules in $\cGr$}

From now on we omit congruence relation symbols for c-groups $A$ and $B$ if no confusion arise.

\begin{defn}
	An (left) \textit{action} of a c-group $B$ on a c-group $A$ is a function $B\times A\rightarrow A$ denoted by $(b,a)\mapsto b\cdot a$ which satisfies the following conditions
	\begin{enumerate}[label=(\roman{*}), leftmargin=1cm]
		\item $b\cdot(a+a_1)\sim (b\cdot a)+(b\cdot a_1)$,
		\item $(b+b_1)\cdot a\sim b\cdot(b_1\cdot a)$,
		\item $0\cdot a\sim a$,
		\item If $a\sim a_1$ and $b\sim b_1$ then $b\cdot a \sim b_1\cdot a_1$,
	\end{enumerate}
	for $a,a_1\in A$ and $b,b_1\in B$.
\end{defn}

A semi-direct product is defined as follow: Let $A, B \in \cGr$ and suppose that $ B $ acts on $ A $ satisfying the conditions (i)--(iv). Then the product $B \times A$ in  $\cGr$ becomes a c-group with the operation $(b',a')+(b,a)=(b'+b, a'+b'\cdot a)$ for any $b, b'\in B, a, a'\in A$ where the congruence relation is the product relation, i.e. $(b,a)\sim (b',a')$ if and only if $b\sim b'$ and $a\sim a'$. Here, $(0,0)$ is a zero element in $B \times A$ and $(-b, -b\cdot (- a))$ is the opposite element of the pair $(b,a)\in B \times A$. This c-group is called the semi-direct product $B \ltimes A$ in $\cGr$.

For any morphism $f:D\rightarrow D'$ in $\cGr$, $f$ is called an isomorphism up to congruence relation or briefly a c-isomorphism if there is a morphism $f':D'\rightarrow D$, such that $ff'\sim 1_{D'}$ and $f'f\sim 1_D$. This kind of an isomorphism is denoted by $\tilde{\approx}$, i.e. $f\colon D \mathrel{\tilde{\approx}} D'$.

For a semi-direct product object $B \ltimes A$ in $\cGr$, we have a natural projection $p': B \ltimes A \rightarrow B$. In this case, there is a c-isomorphism $\cKer p' \mathrel{\tilde{\approx}} A$ which need not to be an isomorphism as in the case of groups.

From a categorical group $\mathsf{C}=(C_0, C_1, d_0, d_1, i, m)$, we obtain a split extension
\begin{equation} \label{extension}
	\xymatrix{0\ar[r]&\cKer d_0\ar[r]^-{j}&C_1\ar[r]_-{d_0}&C_0\ar[r] \ar@/_/[l]_-{i}&0,}
\end{equation}
where $i\colon C_0\rightarrow C_1$ is a section of $d_0$. Thus we define an action of $C_0$ on $\cKer d_0$ by
\[ \begin{array}{rcl}
	C_0\times \cKer d_0& \longrightarrow & \cKer d_0,\\
	(r,c)              & \longmapsto     & r\cdot c=i(r)+(j(c)-i(r)).
\end{array}\]

\begin{prop}
	The action of $C_0$ on $\cKer d_0$ satisfies the conditions for an action in $\cGr$.
\end{prop}

\begin{defn}
	Let $G$ and $H$ be two c-groups, let $\partial\colon G\rightarrow H$ be a morphism of c-groups and let $H$ act on $G$. We call $(G,H,\partial)$ a \emph{c-crossed module} if the following conditions are satisfied:
	\begin{enumerate}[label=(\roman{*}), leftmargin=1cm]
		\item $\partial(b\cdot a)= b+(\partial(a)-b)$,
		\item $\partial(a)\cdot a_1\sim a+(a_1-a)$,
	\end{enumerate}
	for $a,a_1\in G$ and $b\in H$.
\end{defn}

Let $(G,H,\partial)$ and $(G',H',\partial')$ be two c-crossed modules. A \emph{c-crossed module morphism} is a pair of morphisms $\langle f,g\rangle\colon(G,H,\partial)\rightarrow(G',H',\partial')$ such that $g\partial=\partial'f$ and $f(b\cdot a)= g(b)\cdot f(a)$, for all $b\in H$ and $a\in G$, where $f$ and $g$ are morphisms of c-groups. c-crossed modules and morphisms of c-crossed modules form a category.

Let $H$ be a normal c-subgroup of a c-group $G$. One can easily see that, in general, we do not have a usual action by conjugation of $G$ on $H$. However, we have a similar situation as given below.

\begin{lem}
	If $H$ is a perfect normal c-subgroup of a c-group $G$, then we have an action of $G$ on $H$ in the category $\cGr$ and the inclusion morphism defines a c-crossed module.
\end{lem}

\begin{defn}\label{connected cr}
	A c-crossed module $(G,H,\partial)$ is called \textit{connected} if $G$ is a connected c-group.
\end{defn}

\begin{prop}
	For a categorical group $\mathsf{C}=(C_0, C_1, d_0, d_1, i, m)$, $(\cKer d_0,C_0,d)$ is a connected c-crossed module where $d={d_1}|_{\cKer d_0}$.
\end{prop}

In \cite{Dat2020} we introduced another c-group, $\textbf{Star}_{\mathsf{C}} 0$,  from a categorical group $\mathsf{C}=(C_0, C_1, d_0, d_1, i, m)$ apart from $\cKer d_0$. $\textbf{Star}_{\mathsf{C}} 0=\{f\in C_1\mid  d_0(f)=0\}$ with the addition operation given by $f+f'=(f+f')\gamma$, where $f+f':0+0\rightarrow d_1(f)+d_1(f')$ is a sum in $C_1$ and $\gamma : 0 \rightarrow 0+0$ is the unique special isomorphism in $C_1$. The congruence relation on $\textbf{Star}_{\mathsf{C}} 0$ is induced by the relation on $C_1$, which is the relation of being isomorphic in $C_1$. $C_0$ is also a c-group where the congruence relation is given by isomorphisms between the objects. Moreover, there is an action of $C_0$ on $\textbf{Star}_{\mathsf{C}} 0$ given by $r\cdot c=(i(r)+(c-i(r)))\gamma$, for any $r\in C_0, c\in \textbf{Star}_{\mathsf{C}} 0$, where $\gamma\colon 0\approx r+(0-r)$ is a special isomorphism.

\begin{defn}\label{strict cr}
	A c-crossed module $(G,H,\partial)$ is called \textit{strict} if,
	\begin{enumerate}[label=(\roman{*}), leftmargin=1cm]
		\item $\partial(b\cdot a)= b+(\partial(a)-b)$,
		\item $\partial(a)\cdot a_1= a+(a_1-a)$,
	\end{enumerate}
	for $a,a_1\in G$ and $b\in H$.
\end{defn}

\begin{defn}\label{weak special}
	In a c-crossed module $(G,H,\partial)$ a congruence $g\sim g'$ in $G$ is called a \textit{weak special congruence} if $\partial (g)\sim\partial(g')$ is a special congruence in $H$.
\end{defn}

For a c-crossed module $(G,H,\partial)$ every special congruence in $G$ is a weak special congruence since the morphism $\partial$ carries any special congruence to a special congruence between pairs of elements.

\begin{defn}\label{special}
	A c-crossed module $(G,H,\partial)$ is called \emph{special} if for any congruence $\gamma :\partial(c)\sim r$,  there exists $c'\sim c$, such that $\partial(c') =r$, where $c, c'\in G$ and $r\in H$. If $\gamma $ is a special congruence, then  $c'$ is the unique element in $G$ with this property.
\end{defn}
It is easy to show, that if $\gamma $ is a special congruence then the above stated property is equivalent to the following one: there exists a unique element $c'\sim c$ with $\partial(c') =r$ and $c'$ is weak special congruent to $c$.

A c-crossed module is called as a \textit{cssc-crossed module} if it is connected, strict, and special c-crossed module. This kind of crossed module is exactly that we were looking for.

\begin{thm}\label{theo:1}
	For a categorical group $\mathsf{C}=(C_0, C_1, d_0, d_1, i, m)$ the triple $(\textbf{\em Star}_{\mathsf{C}} 0, C_0, d)$ is a cssc-crossed module, where $d={d_1}|_{\textbf{\em Star}_{\mathsf{C}}0}$.
\end{thm}

In \cite{Dat2024}, we constructed a categorical group $\mathbb{G}=(\mathbb{G}_0,\mathbb{G}_1,d_0,d_1,i,m)$ from a cssc-crossed module $\partial: M \rightarrow N$ as follows: Let $N$ be the set of objects $\mathbb{G}_0$ of the desired categorical group, i.e. $\mathbb{G}_0=N$. From the definition of a c-group there exist special congruences between certain elements in $N$. We denote these congruences by $\alpha, \beta,...$. If $\alpha:r \sim  r' $ is a special congruence in $N$, then we write an arrow $\alpha: r \rightarrow  r' $ with $\mathsf{dom}(\alpha)=r$, $\mathsf{codom}(\alpha)=r'$. Since for any special congruence $\alpha: r \sim r' $, there exists the special congruence ${{\alpha}^{-1}}: r' \sim  r $ and the corresponding arrows in $N \times M$ are called as special isomorphisms. Now consider the set $N \bar{\times} M$ with elements of the type $\alpha, \beta,...,\beta (r,c) \alpha $, where $\alpha, \beta,...$ are all special isomorphisms in  $N \times M$, $(r,c) \in N \times M$, and in $\beta (r,c) \alpha$ we have $\mathsf{codom}(\alpha)=r$, $\mathsf{dom}(\beta)=\partial(c)+r$. It is easy to see that any element $(r,c) \in N \times M$ is an element in $N \bar{\times} M$, where $\alpha=1_r: r \xrightarrow{=} r$ and $\beta=1_{\partial(c)+r}: \partial(c)+r \xrightarrow{=} \partial(c)+r$ are identity arrows.

\begin{lem}\label{lem1}
	Let $\partial: M \rightarrow N$ be a c-crossed module. Then we have a special congruence $\partial(0)\sim 0$ in $N$.
\end{lem}

\begin{lem}\label{lem2}
	Let $\partial: M \rightarrow N$ be a cssc-crossed module and $\alpha :r \sim  r' $ a special congruence  in $N$. Then there exists a unique element $c \in M$ with the properties that $c$ is weakly special congruent to $0$ and   $\partial (c)=r'-r$.
\end{lem}

\begin{cor}\label{cor1}
	If there is a special congruence $\alpha :r \sim  r' $ for $r,r' \in N$, then there exists a unique element $c \in M$ with $\partial(c)=r'-r$, a weak congruence $c\sim  0 $ and an arrow $\varepsilon(r,c): r \rightarrow r'$ in $N \bar{\times} M$, where $\varepsilon: (r'-r)+r \rightarrow r'$ is a special isomorphism corresponding to the special congruence $(r'-r)+r\sim r'$ in $N$ and $ (r,c): r \rightarrow (r'-r)+r$ is an arrow in $N \bar{\times} M$.
\end{cor}

In order to define arrows in the desired categorical group, in \cite{Dat2024} we defined two identifications in $N \bar{\times} M$ as follows:

\begin{description}
\item[Identification I in $N \bar{\times} M$:] In what follows we will identify  any special isomorphism $\alpha :r \sim  r' $ in $N \bar{\times} M$ with the arrow $\varepsilon(r,c)$  defined uniquely; i.e. $c$ is the element defined uniquely and weakly special congruent to $0$ in $M$ with $\partial (c)=r'-r$ and $\varepsilon$ is a  unique special isomorphism defined in Corollary \ref{cor1}. We will use the notation $\alpha \equiv \varepsilon(r,c)$.
\item[Identification II in $N \bar{\times} M$:] If $\varphi:c\sim c'$ is a weak special congruence in $M$, then $\beta (r,c)\alpha\equiv\beta' (r,c')\alpha$, for any $r\in N$, where $\beta=\beta'(\partial(\varphi)+1_r)$, for any special congruence $\beta'$ with $\mathsf{dom}(\beta')=\partial(c')+r$, and any $\alpha$ with $\mathsf{codom}(\alpha)=r$
 \begin{equation}\label{diag:31}
  \begin{gathered}
	\xymatrix@R=15mm@C=15mm{
		&  & \partial(c)+r \ar[dr]^-{\beta}\ar[dd]^-{\approx}_-{\partial(\varphi)+1_r} &  \\
		r' \ar[r]^-{\alpha} & r \ar[ur]^-{(r,c)} \ar[dr]_-{(r,c')} &  & r'' \\
		&  & \partial(c')+r \ar[ur]_-{\beta'} &
		}
  \end{gathered}
 \end{equation}
\end{description}

In the diagram, in particular,  if   $\beta'=1_{\partial(c')+r}$ then $\beta=\partial(\varphi)+1_r$ and we obtain $(r,c')\alpha=(\partial(\varphi)+1_r)(r,c)\alpha$. If in addition $\alpha=1_r$, then we have $(r,c')=(\partial(\varphi)+1_r)(r,c)$.

Obviously, if $\alpha\equiv\varepsilon(r,c)$ by \textbf{Identification I} and $\varepsilon(r,c)\equiv\varepsilon'(r,c')$ by \textbf{Identification II}, then $\alpha\equiv\varepsilon'(r,c')$; also if $\beta(r,c)\alpha\equiv\beta'(r,c')\alpha\equiv\beta''(r,c'')\alpha$ by \textbf{Identification II}, then $\beta(r,c)\equiv\beta''(r,c'')\alpha$.

We denoted the resulting set obtained by \textbf{Identifications I} and \textbf{II} by $\mathbb{G}_{1}$; the elements of $\mathbb{G}_{1}$ is called arrows of the categorical group $\mathbb{G}$.

Compositions of arrows in $\mathbb{G}_1$ were defined in the following ways:

\begin{defn} \label{def1} \begin{enumerate}[label=(\textit{\alph{*}}), leftmargin=1cm]
		\item\label{def1:a} Let $\alpha$ and $ \beta$  be two composable special isomorphisms in $\mathbb{G}_{1}$. The composition  $ \beta \circ \alpha$ is defined to be $ \beta\alpha$ which is the special isomorphism corresponding to the composition of special congruences  $\alpha$ and $ \beta$ in $N$ which is a special congruence as well.
		\item\label{def1:b} If $\beta$ and $ \gamma$  are composable special isomorphisms in $\mathbb{G}_{1}$, the composition $\gamma \circ(\beta (r,c)\alpha)$ is defined as $ (\gamma\beta)(r,c)\alpha$. In analogous way  $(\beta (r,c)\alpha) \circ \delta=\beta (r,c)(\alpha\delta)$, where $\delta$ and $\alpha$ are composable.
		From this definition it follows that
		\begin{equation}\label{eq:1}
			\gamma \circ ((r,c)\alpha)=\gamma (r,c)\alpha
		\end{equation}
		\begin{equation}\label{eq:23}
			(\beta (r,c)) \circ \delta=\beta (r,c) \delta
		\end{equation}
		\item\label{def1:c} Let $\beta (r,c)\alpha$ and $\beta' (r',c')\alpha'\in \mathbb{G}_{1}$. If $\beta$ and $\alpha'$ are composable, we define
		\begin{equation*}
			(\beta' (r',c')\alpha') \circ (\beta (r,c)\alpha)=\theta (r,c'+c)\alpha,
		\end{equation*}
		where $\theta$ is the composition of the following special isomorphisms
		\begin{equation*}
			(\partial(c')+\partial(c))+r\approx \partial(c')+(\partial(c)+r)\overset{{{1}_{\partial(c')}}+\beta }{\mathop{\approx }}\,\partial(c')+{{r}_{2}}\overset{{{1}_{\partial(c')}}+\alpha }{\mathop{\approx }}\,\partial(c')+r'\overset{\beta }{\mathop{\approx }}\,r''
		\end{equation*}
		the picture is
		\begin{equation*}
			{{r}_{1}}\xrightarrow{\alpha }r\xrightarrow{(r,c)}\partial(c)+r\xrightarrow{\beta }{{r}_{2}}\xrightarrow{\alpha '}r'\xrightarrow{(r',c')}\partial(c')+r'\xrightarrow{\beta '}r''.
		\end{equation*}
	\end{enumerate}
\end{defn}

It is proved that definitions of compositions do not depend on the choice of representatives \cite[Propositions 4.2-4.5]{Dat2024}. Also, in \cite[Proposition 4.6]{Dat2024} we proved that composition of arrows is associative.

\begin{prop}\label{prop:46}
	Let $(r,c),(r',c') \in \mathbb{G}_{1}$ and $r \sim r'$ be a congruence in $N$. Then there exist arrows $\varphi: r \rightarrow r'$ and $\theta:\partial(c)+r \rightarrow \partial(c')+r'$ for which the following diagram commutes
	\begin{equation*}
		\xymatrix@R=15mm@C=15mm{
			r \ar[r]^-{(r,c)} \ar[d]_-{\varphi} & \partial(c)+r \ar[d]^-{\theta} \\
			r' \ar[r]_-{(r',c')} & \partial(c')+r'
		}
	\end{equation*}
	%		\[
	%		\begin{tikzpicture}
		%			\node (s) {$r$};
		%			\node (xy) [below=2 of s] {$r' $};
		%			\node (x) [right=3 of xy] {$\partial(c')+r'$};
		%			\node (y) [right=3 of s] {$(\partial(c)+r)$};
		%			\draw[->] (s) to node [sloped, above] {$(r,c)$} (y);
		%			\draw[->] (s) to node [anchor=east]{$\psi$} (xy);
		%			\draw[->] (xy) to node [below]{$(r',c')$} (x);
		%			\draw[->] (y) to node [below][anchor=west] {$\theta$} (x);
		%		\end{tikzpicture}
	%		\]
\end{prop}

\begin{cor}\label{cor:2}
	Let $\beta (r,c)\alpha$ and $ \beta' (r',c')\alpha' \in \mathbb{G}_{1}$, where $\mathsf{dom}(\alpha)=\mathsf{dom}(\alpha')$ and $\mathsf{codom}(\beta)=\mathsf{codom}(\beta')$. If $r\sim r'$ is a congruence and $c\sim c'$ is a weak special congruence, then $\beta (r,c)\alpha=\beta' (r',c')\alpha'$.
\end{cor}

For the identity arrows for any object $r \in\mathbb{G}_0$, we defined an arrow $\varepsilon (r,0)$, where $(r,0): r \rightarrow \partial(0)+r$ is an arrow in $\mathbb{G}_1$ and $\varepsilon:\partial(0)+r \rightarrow r $ is a special isomorphism by means of the special isomorphism $\partial(0)\approx 0$. We showed that $\varepsilon (r,0):r \rightarrow r $ is an identity arrow. Moreover, for the inverse of any arrow $\beta (r,c)\alpha \in \mathbb{G}_1$ define an arrow $\alpha^{-1} (r,c)^{-1}\beta^{-1}$, where $(r,c)^{-1}=\varphi(\partial(c)+r,-c)$, and $\varphi:\partial(-c)+(\partial(c)+r)\rightarrow r$ is a special isomorphism. In this case,

\begin{equation*}
	(\beta (r,c)\alpha)(\alpha^{-1}\varphi(\partial(c)+r,-c)\beta^{-1})=1_{r''},
\end{equation*}
\begin{equation*}
	(\alpha^{-1}\varphi(\partial(c)+r,-c)\beta^{-1})(\beta (r,c)\alpha)=1_{r'}	
\end{equation*}
where $\alpha\colon r'\rightarrow r$ and $\beta\colon \partial(c)+r\rightarrow r''$.

It is also proved that the inverse of an arrow does not depend on the choice of representatives \cite[Proposition 5.3]{Dat2024}.

\begin{prop}
	If $\varphi:c\sim c'$ is a weak special congruence, then according to diagram \eqref{diag:31} for any $r \in \mathbb{G}_0$ we have
	\begin{equation*}
		(\beta' (r,c')\alpha)^{-1}=(\beta (r,c)\alpha)^{-1}.
	\end{equation*}
\end{prop}

Addition operation for any arrows $\beta_1(r_1,c_1)\alpha_1$ and $ \beta_2(r_2,c_2)\alpha_2\in\mathbb{G}_1$ is defined by

\begin{equation*}
	\beta_1(r_1,c_1)\alpha_1+\beta_2(r_2,c_2)\alpha_2=(\beta_1+\beta_2)\theta(r_1+r_2,c_1+r_1\cdot c_2)(\alpha_1+\alpha_2),
\end{equation*}
where
\begin{equation*}
	\theta:\partial(c_1)+(r_1+(\partial(c_2)-r_1))+r_1+r_2\longrightarrow(\partial(c_1)+r_1)+(\partial(c_2)+r_2)	
\end{equation*}
is the special isomorphism.

Let $\alpha\colon r\rightarrow r'$ be a special isomorphism in $\mathbb{G}_1$. By \textbf{Identification I}, we have $\alpha\equiv\varepsilon(r,c)$, where $c\sim 0$ is a weak special congruence with $\partial(c)=r'-r$ and $\varepsilon\colon (r'-r)+r\rightarrow r'$ is the special isomorphism. For any $\gamma(r',c')\beta\in\mathbb{G}_1$ the sum $\alpha+\gamma(r',c')\beta$  is defined to be $\varepsilon(r,c)+\gamma(r',c')\beta$. Similarly, $\gamma(r',c')\beta+\alpha$ is defined to be $\gamma(r',c')\beta+\varepsilon(r,c)$.

It is also proved that the addition operation in $\mathbb{G}_1$ does not depend on the choice of representatives \cite[Propositions 6.1,6.2]{Dat2024} and is associative up to isomorphism \cite[Proposition 6.3]{Dat2024}.

The sum of objects in $\mathbb{G}_0$ is defined as the sum in c-group $N$, which is associative up to special congruence, and therefore it is associative in $\mathbb{G}_0$ up to special isomorphism in $\mathbb{G}_1$. A zero object in the category $\mathbb{G}$ with objects $\mathbb{G}_0=N$ is a zero element 0 in $N$ as in c-group, i.e. we have special isomorphisms $0+r\approx r\approx r+0$ which are defined by special congruences in $N$. A zero element for the addition operation in $\mathbb{G}_1$ we defined by $\sigma(0,0)$, where $\sigma$ is the special isomorphism defined by the special congruence $\sigma\colon \partial(0)+0\sim 0$. According to the definition of identity arrows for any object $r\in\mathbb{G}_0$, we have $\sigma(0,0)=1_0\colon 0\rightarrow 0$.

For any arrow $\beta(r,c)\alpha\in\mathbb{G}_1$ we have an isomorphism
\begin{equation*}
	\beta(r,c)\alpha+\sigma(0,0)\approx  \beta(r,c)\alpha \approx \sigma(0,0)+\beta(r,c)\alpha,
\end{equation*}
and we defined an opposite arrow $-(\beta(r,c)\alpha)=(-\beta)\varepsilon(-r,(-r)\cdot (-c))(-\alpha)$, where $\alpha\colon r'\rightarrow r$ and $\beta\colon \partial(c)+r\rightarrow r''$ are special isomorphisms. $-\alpha\colon -r'\rightarrow -r$ and $-\beta\colon -(\partial(c)+r)\rightarrow -r''$ are the opposite special isomorphisms  corresponding to the opposite special congruences $-\alpha\colon -r'\sim  -r$ and $-\beta\colon -(\partial(c)+r)\sim -r''$, respectively; and $\varepsilon\colon ((-r)+(\partial(-c)-(-r)))+(-r)\rightarrow -(\partial(c)+r)$ is an obvious special isomorphism.

\begin{prop}\label{prop:oppoarrow}
For any $\beta(r,c)\alpha\in\mathbb{G}_1$ we have isomorphisms
 \begin{equation*}
	\beta(r,c)\alpha+(-(\beta(r,c)\alpha))\approx \sigma(0,0)\approx -(\beta(r,c)\alpha)+\beta(r,c)\alpha
 \end{equation*}
where $\sigma\colon \partial(0)+0\rightarrow 0$ is the special isomorphism.
\end{prop}

An opposite arrow does not depend on the choice of representatives \cite[Proposition 6.6]{Dat2024}.

\begin{thm}\label{theo:2}
The category $\mathbb{G}=(\mathbb{G}_0,\mathbb{G}_1,d_0,d_1,i,m)$ with the defined addition operation in $\mathbb{G}$ is a categorical group.
\end{thm}

Denote by $\Catgr$ and $\cCrsM$ the categories of categorical groups and  cssc-crossed modules respectively.

\section{Functoriality of the correspondence $\Catgr\longrightarrow\cCrsM$}

For any categorical group $\mathsf{C}=(C_0,C_1,d_0,d_1,i,m)$ according to Theorem \ref{theo:1} we define a correspondence $\mathbb{L}_0 \colon \Catgr \longrightarrow \cCrsM$ between the objects of these categories by $\mathbb{L}_0(\mathsf{C})=(\Str,C_0,d)$, where $d=d_1|_{\Str}$.

Let $\mathsf{C'}=(C'_0, C'_1, d'_0, d'_1, i', m')$ be another categorical group and $F=(F_0,F_1)\colon \mathsf{C}\longrightarrow \mathsf{C}'$ be a functor between them. We define a map $\varphi_0\colon C_0\longrightarrow C'_0 $ as $F_0$ and define a map $\varphi_1\colon  \Str\longrightarrow \Strp$ as $\varphi_1=F_1\mid_{\Str}$. Therefore by definition $\mathbb{L}_1(F)=(\varphi_0,\varphi_1)$.  It follows that $\varphi_0$ and $\varphi_1$ are morphisms of c-groups and the diagram
\begin{equation*}
		\xymatrix@R=15mm@C=15mm{
		\Str \ar[r]^-{ \partial} \ar[d]_-{\varphi_1} & C_0 \ar[d]^-{\varphi_0} \\
		\Strp  \ar[r]_-{\partial'} & C'_0
		}
\end{equation*}
is commutative.

\begin{prop}
$\left(\varphi_0, \varphi_1\right)$ is a morphism of c-crossed modules.
\end{prop}

\begin{proof}
We have to check, that $\varphi_1(r\cdot f)=\varphi_0(r)\cdot \varphi_1(f)$, for any $r\in C_0$ and $f\colon 0 \rightarrow r_1$ from $\Str$. By the definition of action of $C_0$ on $\Str$ we have $r\cdot f=(i(r)+(f-i(r)))\gamma$, where $\gamma\colon  0_{\mathsf{C}} \xrightarrow{~\approx~} r+(0-r)$ is a special isomorphism in $\mathsf{C}$.

  Similarly, $\varphi_0(r)\cdot\varphi_1(f)=(i'\varphi_0(r)+(\varphi_1(f)-i'\varphi_0(r)))\gamma'$, where $\gamma'\colon  0_{\mathsf{C}'} \xrightarrow{~\approx~~} \varphi_0(r)+(\varphi_0(0)-\varphi_0(r))$ is a special isomorphism in $\mathsf{C'}$. We have to check, that
  \begin{equation}\label{eq:31}
		\varphi_1((i(r)+(f-i(r)))\gamma)=(i'\varphi_0(r)+(\varphi_1(f)-i'\varphi_0(r)))\gamma'.
		\end{equation}
  The left side is equal to $(i'\varphi_0(r)+(\varphi_1(f)-i'\varphi_0(r)))\varphi_1(\gamma)$, where $\varphi_1(\gamma)=\gamma'$, since $\C'$ is coherent, which gives equality \eqref{eq:31}.
\end{proof}

Denote the constructed correspondence by $\mathbb{L}=(\mathbb{L}_0,\mathbb{L}_1)\colon  \Catgr \longrightarrow \cCrsM$. It is easy to show that this correspondence carries $1_{\C}\colon\C\rightarrow \C$ to $1_{\Str\rightarrow C_0} \colon  (\Str\rightarrow C_0)\longrightarrow (\Str\rightarrow C_0)$ and for the composition $\mathsf{C}\xrightarrow{~F~}\mathsf{C'} \xrightarrow{~P~}\mathsf{C''}$ we have $ PF\longmapsto \mathbb{L}(P)\cdot\mathbb{L}(F)$. Therefore we proved that $\mathbb{L}$ is a functor.
	
\section{Functoriality of the correspondence $\cCrsM\longrightarrow\Catgr$}

For any cssc-crossed module $\partial\colon M\rightarrow N$ according to Theorem \ref{theo:2} we define a correspondence \linebreak $\mathbb{T}_0 \colon \cCrsM\longrightarrow\Catgr$ between the objects of these categories by $\mathbb{T}_0(\partial\colon M\rightarrow N)=\mathbb{G}=({\mathbb{G}}_{0},{\mathbb{G}}_{1},d_0,d_1,i,m)$.
 
Let $\partial'\colon  M'\rightarrow N'$  be another cssc-crossed module and $ \varphi=(\varphi_0,\varphi_1)$ be a morphism between them
\begin{equation*}
		\xymatrix@R=15mm@C=15mm{
		M \ar[r]^-{ \partial} \ar[d]_-{\varphi_1} & N \ar[d]^-{\varphi_0} \\
		M' \ar[r]_-{\partial'} & N'.
		}
\end{equation*}

We have to construct a pair of maps $\mathbb{T}_{1,0}\colon\mathbb{G}_0\rightarrow \mathbb{G}'_0$ and $\mathbb{T}_{1,1}\colon \mathbb{G}_1\rightarrow \mathbb{G}'_1$ and show that they define a functor $\mathbb{T}_0\left(M\xrightarrow{~\partial~} N\right)\xlongrightarrow{~\mathbb{T}_1(\varphi)~}\mathbb{T}_0\left(M'\xrightarrow{~\partial'~} N'\right) $ between the categorical groups, where $\mathbb{T}_0(\partial\colon  M'\rightarrow N')=\mathbb{G}'=(\mathbb{G}'_0,\mathbb{G}'_1,d'_0,d'_1,i',m')$. We define ${\mathbb{T}}_{1,0}(\varphi)=\varphi_0\colon \mathbb{G}_0\rightarrow \mathbb{G}'_0$ and ${\mathbb{T}}_{1,1}(\varphi)\colon \mathbb{G}_1\rightarrow \mathbb{G}'_1$ is defined in a natural way $\beta(r,c)\alpha\longmapsto\varphi_1(\beta)(\varphi_0(r),\varphi_1(c))\varphi_1(\alpha)$ for any $\beta(r,c)\alpha\in \mathbb{G}_1$. Therefore we defined ${\mathbb{T}}_{1}(\varphi)=({\mathbb{T}}_{1.0}(\varphi),{\mathbb{T}}_{1,1})(\varphi)$.

Below for simplicity we will not use indices for the defined functors. We will write $\mathbb{T}(\varphi_0)$ and  $\mathbb{T}(\varphi_1)$ instead of
${\mathbb{T}}_{1,0}(\varphi_0,\varphi_1)$ and  ${\mathbb{T}}_{1,1}(\varphi_0,\varphi_1)$ respectively.

\begin{prop}\label{prop:41}
  $\left(\mathbb{T}(\varphi_0), \mathbb{T}(\varphi_1)\right)$ is a morphism of categorical groups.
\end{prop}

\begin{proof}
It is obvious that $\mathbb{T}(\varphi_0)$ preserves sums, the zero element and opposite elements. We have to prove that
\begin{alignat}{2}
\mathbb{T}(\varphi_1)(\beta(r,c)\alpha+\beta'(r',c')\alpha')&=\mathbb{T}(\varphi_1)(\beta(r,c)\alpha)+\mathbb{T}(\varphi_1)(\beta'(r',c')\alpha') \label{eq:2}\\
\mathbb{T}(\varphi_1)(\epsilon(0,0))&=(\epsilon'(0',0')) \label{eq:3} \\
\mathbb{T}(\psi\varphi)&=\mathbb{T}(\psi)\mathbb{T}(\varphi) \label{eq:4}
\end{alignat}
%\begin{equation}\label{eq:2}\mathbb{T}(\varphi_1)(\beta(r,c)\alpha+\beta'(r',c')\alpha')=\mathbb{T}(\varphi_1)(\beta(r,c)\alpha)+\mathbb{T}(\varphi_1)(\beta'(r',c')\alpha')
%\end{equation}
%\begin{equation}\label{eq:3}
%\mathbb{T}(\varphi_1)(\epsilon(0,0))=(\epsilon'(0',0'))
%\end{equation}
%\begin{equation}\label{eq:4}
%\mathbb{T}(\psi\varphi)=\mathbb{T}(\psi)\mathbb{T}(\varphi)
%\end{equation}
for any $\beta(r,c)\alpha$ and $\beta'(r',c')\alpha'\in \mathbb{G}_1 $, for $\epsilon(0,0)$ the zero element in $\mathbb{G}_1$, where $\epsilon'(0',0')$ is a zero element in $\mathbb{G'}_1$;  $\psi=(\psi_1,\psi_0)\colon\left(M'\xrightarrow{~\partial'~} N'\right)\longrightarrow \left(M''\xrightarrow{~\partial''~} N''\right) $ and  $\varphi=(\varphi_1,\varphi_0)\colon\left(M\xrightarrow{~\partial~} N\right)\longrightarrow \left(M'\xrightarrow{~\partial'~} N'\right) $ are any composable arrows between the cssc-crossed modules.

\eqref{eq:3} follows from the definition of $\mathbb{T}(\varphi_1)$ and properties of $\varphi_0$ and $\varphi_1$.

For \eqref{eq:2} we compute the both sides of the equality.
\begin{equation}\label{eq:5}
\begin{alignedat}{2}
\mathbb{T}(\varphi_1)(\beta(r,c)\alpha+\beta'(r',c')\alpha')&=\mathbb{T}(\varphi_1)\left((\beta+\beta')\theta(r+r',c+r\cdot c')(\alpha+\alpha')\right)\\
&=\varphi_1(\beta+\beta')\varphi_1(\theta)\left(\varphi_0(r+r'),\varphi_1(c+r\cdot c')\right)\varphi_1(\alpha+\alpha'),
\end{alignedat}
\end{equation}
where $\theta$ is the special congruence defined by the definition of sum in $\mathbb{G}_1$.
\begin{equation}\label{eq:6}
\begin{alignedat}{2}
\mathbb{T}(\varphi_1)(\beta(r,c)\alpha) + \mathbb{T}(\varphi_2)(\beta'(r',c')\alpha') & = \varphi_1(\beta)(\varphi_0(r),\varphi_1(c))\varphi_1(\alpha) +  \varphi_1(\beta')(\varphi_0(r'),\varphi_1(c'))\varphi_1(\alpha') \\
& =(\varphi_1(\beta)+\varphi_1(\beta'))\theta_1(\varphi_0(r)+\varphi_0(r'), \\
& ~~~~\varphi_1(c)+\varphi_0(r)\cdot\varphi_1(c'))\cdot(\varphi_1(\alpha) +  \varphi_1(\alpha')  ),
\end{alignedat}
\end{equation}
where $\theta_1$ is the appropriate special congruence defined according to the definition of a sum in $\mathbb{G}_1$.

We have
\begin{alignat*}{2}
\varphi_1(\beta+\beta') &= \varphi_1(\beta)+\varphi_1(\beta'),\\
\varphi_1(\alpha+\alpha') &= \varphi_1(\alpha)+\varphi_1(\alpha')
\end{alignat*}
and
\begin{equation*}
\varphi_1(c+r\cdot c')=\varphi_1(c)+\varphi_0(r)\cdot \varphi_1(c').
\end{equation*}

From the uniqueness of special congruences we have $\varphi_1(\theta)=\theta_1$. From these equalities follows that expressions \eqref{eq:5} and \eqref{eq:6} are equal, which is the required equality \eqref{eq:2}.

Now we shall show \eqref{eq:4}. The picture is the following
\begin{equation*}
\psi\varphi\colon \begin{cases*}
\xymatrix@R=15mm@C=15mm{
M \ar[r]^-{ \partial} \ar[d]|-{\varphi_1} \ar@/_1.7pc/[dd]_-{\psi_1\varphi_1} & N \ar[d]|-{\varphi_0} \ar@/^1.7pc/[dd]^-{\psi_0\varphi_0}\\
M' \ar[r]^-{\partial'} \ar[d]|-{\psi_1}& N'\ar[d]|-{\psi_0}\\
M''\ar[r]^-{\partial''}& N'',}
\end{cases*}
\end{equation*}

\begin{equation*}
\mathbb{T}(\psi\varphi)\colon\begin{cases}
\xymatrix@C=15mm@R=15mm{
(\mathbb{G}_0, \mathbb{G}_1, d_0, d_1, i, m) \ar[d]^-{\mathbb{T}(\varphi_0,\varphi_1)} \ar@<10ex>@/^1.7pc/[dd]^-{\mathbb{T}\left(\psi\right)\mathbb{T}\left(\varphi\right)} & \\
(\mathbb{G'}_0, \mathbb{G'}_1, d'_0, d'_1, i', m')\ar[d]^-{\mathbb{T}(\psi_0,\psi_1)}& \\
(\mathbb{G''}_0, \mathbb{G''}_1, d''_0, d''_1, i'', m'').
}
\end{cases}
\end{equation*}
%\begin{equation*}
%\xymatrix@C=15mm@R=15mm{
%\mathbb{T}(\psi\varphi)\colon & (\mathbb{G}_0, \mathbb{G}_1, \partial_0, \partial_1, i, m) \ar[d]^-{\mathbb{T}(\varphi_0,\varphi_1)} & \\
%& (\mathbb{G'}_0, \mathbb{G'}_1, \partial'_0, \partial'_1, i', m')\ar[d]^-{\mathbb{T}(\psi_0,\psi_1)}& \\
%& (\mathbb{G''}_0, \mathbb{G''}_1, \partial''_0, \partial''_1, i'', m'').
%}
%\end{equation*}
We have
\begin{alignat*}{2}
\mathbb{T}(\psi_0)\mathbb{T}(\varphi_0)(r)&=\psi_0\varphi_0(r) \\
\mathbb{T}(\psi_0\varphi_0)(r)&=\psi_0\varphi_0(r).
\end{alignat*}
Moreover
\begin{alignat*}{2}
\mathbb{T}(\psi)\mathbb{T}(\varphi)(\beta(r,c)\alpha)&=\mathbb{T}(\psi)(\varphi_1(\beta)(\varphi_0(r),\varphi_1(c)))\varphi_1(\alpha)\\
&=\psi_1\varphi_1(\beta)(\psi_0\varphi_0(r),\psi_1\varphi_1(c))\psi_1\varphi_1(\alpha).
\end{alignat*}
\begin{equation*}
\mathbb{T}(\psi\varphi)(\beta(r,c)\alpha)=\psi_1\varphi_1(\beta)(\psi_0\varphi_0(r),\psi_1\varphi_1(c))\psi_1\varphi_1(\alpha),
\end{equation*}
which proves equality \eqref{eq:4} and Proposition \ref{prop:41}.
\end{proof}

For any cssc-crossed module $\partial:M\rightarrow N$ and the identity morphism $(1_N,1_M)$ we have $\mathbb{T}(1_N,1_M)=1_{\mathbb{T}\left(M\xrightarrow{~\partial~} N\right)}$, where $1_{\mathbb{T}\left(M\xrightarrow{~\partial~} N\right)}$ is an identity functor $\mathbb{T}\left(M\xrightarrow{~\partial~} N\right)\xrightarrow{~=~} \mathbb{T}\left(M\xrightarrow{~\partial~} N\right) $. This note completes the proof of functoriality of the correspondence $\cCrsM\longrightarrow \Catgr$.

\section{The isomorphism $\mathbb{TL}\approx 1_{\Catgr}$}

Let $\C=(C_0, C_1, d_0, d_1, i, m)$ be a categorical group. As we know \cite{Dat2024} (see section 3) $\mathbb{L}(\C)$ is a crossed module $\partial\colon\Str\rightarrow C_0$, where $\partial=d_1\mid_{\Str}$ and $\mathbb{T}\left(\Str\xrightarrow{~\partial~} C_0\right)=\mathsf{\overline{C}}=(\mathbb{G}_0, \mathbb{G}_1, \overline{d}_0, \overline{d}_1, \overline{i}, \overline{m})$. We will construct functors between the categorical groups
\begin{equation*}
\xymatrix@R=15mm@C=15mm{
\C \ar@<.7ex>[r]^-{P} & {\overline{\C}} \ar@<.3ex>[l]^-{F}}
\end{equation*}
and prove that $FP=1_{\C}$ and $PF=1_{\overline{\C}}$.

Obviously, $P_0\colon  C_0\rightarrow \mathbb{G}_0$ is the identity map, which is a morphism of c-groups. For any $f\in C_1$ we define $P_1(f)=\theta(d_0f,(f-id_0f)\alpha)$, where $\alpha$ is a special isomorphism in $C_1$, and $\theta$ is defined as isomorphism $\partial(f-id_0f)\alpha+d_0f=(d_1f-d_0f)+d_0f\xrightarrow{~~\theta~~} d_1f$. For any special isomorphism $\alpha$ in $C_1$, we define $P_1(\alpha)=\alpha$.
\begin{prop}\label{prop:51}
$(P_0,P_1)$ is a morphism of categorical groups.
\end{prop}
\begin{proof}
The proof is based on several lemmas.
\end{proof}

\begin{lem}\label{lem:52}
$P_0$ and $P_1$ preserve the structures in $C_0$ and $C_1$, respectively.
\end{lem}

\begin{proof}
As we have noted already $P_0$ is a morphism of c-groups. Therefore we have only to show that $P_1(f+g)=P_1(f)+P_1(g)$.

$P_1(0)=0$, from which it follows that $P_1(-f)=-P_1(f)$.

We have $P_1(f+g)=\theta_1\left(d_0(f+g),((f+g)-id_0(f+g))\overline{\alpha}\right)$, where $\overline{\alpha}\colon 0\rightarrow d_0(f+g)-d_0(id_0(f+g))$ is an obvious special isomorphism and $\theta_1$ is defined by
\begin{equation*}
d_1((f+g)-id_0(f+g))\overline{\alpha})+d_0(f+g)=(d_1f+d_1g)-(-d_0g-d_0f)+(d_0f+d_0g)\stackrel{~\theta_1~}{\approx} d_1f+d_1g.
\end{equation*}
We have
\begin{alignat*}{2}
P_1(f)+P_1(g)& =\eta_1(d_0f,(f-id_0f)\alpha_1)+\eta_2(d_0g,(g-id_0g)\alpha_2)\\
& =(\eta_1+\eta_2)\xi(d_0f+d_0g,(f-id_0f)\alpha_1 +d_0f\cdot(g-id_0g)\alpha_2)\\
& =(\eta_1+\eta_2)\xi(d_0f+d_0g,(f-id_0f)+(id_0f+(g-id_0g)-id_0f)\alpha_3),
\end{alignat*}
where $\alpha_3$ is defined in the usual way. Now we have only to note, that $d_0(f+g)=d_0f+d_0g$ and $((f+g)-id_0(f+g))\overline{\alpha})$ is special isomorphic to $(f-id_0f)+(id_0f+(g-id_0g)-id_0f)\alpha_3$. Therefore by Corollary \ref{cor:2} follows that $P_1(f+g)=P_1(f)+P_1(g)$.

For the next equality of the lemma we have $$P_1(0)=(0,(0-id_0 0)\overline{\overline{\alpha}}))=\eta(0,0),$$ where $\overline{\overline{\alpha}}$ and $\eta$ are defined in the appropriate way, which proves the lemma.
\end{proof}

\begin{lem}\label{lem:53}
$P_1(i(r))$ is the identity arrow in $\mathsf{\overline{C}}$, for any object $r\in C_0$, and $P_1(fg)=P_1(f)\cdot P_1(g)$.
\end{lem}
\begin{proof}
We have $P_1(i(r))=\theta (d_0i(r),(i(r)-id_0(i(r))\Tilde{\alpha}))=\eta(r,0)$   by Corollary \ref{cor:2} We have $P_1(fg)=\theta(d_0(fg),(fg-id_0(fg))\alpha)=\theta(d_0g,(fg-id_0(g))\alpha)$, where $\theta$ and $\alpha$ are the appropriate special isomorphisms. For the right side of the equality in the lemma we have
\begin{alignat*}{2}
P_1(f)\cdot P_1(g) &= \theta_1(d_0f,(f-id_0f)\alpha_1)\cdot\theta_2(d_0g,(g-id_0g)\alpha_2)\\
&=\overline{\theta}_1(d_0g,((f-id_0f)+(g-id_0g))\overline{\alpha}),
\end{alignat*}
where $\overline{\theta}_1$ and $\overline{\alpha}$ are defined in the appropriate way. By Lemma \ref{lem:compadd} we have the special isomorphism $fg\approx(f-id_0f)+g$, from which we obtain following special isomorphisms
\begin{equation*}
(fg-id_0g)\alpha\approx((f-id_0f)+g)-id_0g\approx(f-id_0f)+(g-id_0g),
\end{equation*}
which by Corollary \ref{cor:2} proves that $P_1(fg)=P_1(f)\cdot P_1(g)$.
\end{proof}

\begin{prf}
By Lemmas \ref{lem:52} and \ref{lem:53} it follows that $(P_0,P_1)$ is a morphism of categorical groups. \hfill \qed
\end{prf}

Define $F_0(r)=r$, for any object $r$ in $\overline{\C}$, $F_1(\eta)=\eta$ for any special isomorphism $\eta$ in $\overline{\C}$, $F_1\left(\beta(r,c)\alpha\right)=\beta\left(c+i(r)\right)\varepsilon\alpha$, for any arrow $\beta(r,c)\alpha$ in $\overline{\C}$, where $\alpha\colon r'\rightarrow r$, $\beta\colon\partial(c)+r\rightarrow r''$ and $\varepsilon\colon r\rightarrow 0+r$ are special congruences in $\Str$.

\begin{lem}\label{lem:54}
$F_0$ and $F_1$ preserve the structures in $\mathbb{G}_0$ and $\mathbb{G}_1$, respectively.
\end{lem}
\begin{proof}
Since $F_0$ is the identity map, we need only to prove that
\begin{equation*}
F_1\left(\beta\left(r,c\right)\alpha+\beta'\left(r',c'\right)\alpha'\right)=F_1\left(\beta\left(r,c\right)\alpha\right)+F_1\left(\beta'\left(r',c'\right)\alpha'\right)
\end{equation*}
and $F_1\left(\xi(0,0)\right)=0$, where $\xi\colon\partial(0)+0\sim 0$ is the special congruence. We compute both sides of the first equality.

 For the left side we have
\begin{alignat*}{2}
F_1\left(\beta\left(r,c\right)\alpha+\beta'\left(r',c'\right)\alpha'\right) & = F_1\left(\left(\beta+\beta'\right)\theta\left(r+r',c+r\cdot c'\right)\left(\alpha+\alpha'\right)\right) \\
& = \left(\beta+\beta'\right)\left(\left(c+r\cdot c'\right)+i\left(r+r'\right)\right)\overline{\varepsilon}\left(\alpha+\alpha'\right),
\end{alignat*}
and
\begin{alignat*}{2}
F_1\left(\beta\left(r,c\right)\alpha\right)+F_1\left(\beta'\left(r',c'\right)\alpha'\right) & = \beta\left(c+i(r)\right)\varepsilon\alpha +  \beta'\left(c'+i(r')\right)\varepsilon'\alpha' \\
& \approx \left(c+i(r)\right)+\left(c'+i(r')\right) \\
& \approx c+\left(i(r)+c'-i(r)\right)+\left(i(r)+i(r')\right) \\
& \approx \left(c+r\cdot c'\right)+i\left(r+r'\right).
\end{alignat*}

Here we have to remember that $\overline{\C}$ is obtained from $\C$ by the correspondences
\begin{equation*}
\C \longmapsto \left(\Str\rightarrow C_0\right) \longmapsto \overline{\C}.
\end{equation*}

From the coherence property of $\C$ we obtain the desired equality. For the next equality, we have
\begin{equation*}
F_1(\xi(0,0)) = \xi\left(1_0+\overline{i}(0)\right)\overline{\varepsilon}=0_{C_{1}}.
\end{equation*}
\end{proof}

\begin{lem}\label{lem:55}
$F_1$ carries identity arrows to the identity arrows and compositions of arrows to the compositions.
\end{lem}
\begin{proof}
We have to prove
\begin{equation*}
F_1\left( \overline{i}(r)\right)  = i(r)
\end{equation*}
and
\begin{equation}\label{eq:51}
F_1 \left( \beta'(r''', c') \alpha' \circ \beta(r, c)\alpha \right) = F_1 \left( \beta'(r''', c')\alpha' \right)\cdot F_1(\beta(r, c)\alpha).
\end{equation}
We have
\begin{equation*}
F_1(\overline{i}(r)) = F_1 \left( \xi(r, 0) \right) = \xi( 1_0 + \overline{i}(r))\varepsilon = i(r)
\end{equation*}
by coherence property of $\C$. Here $ \xi\colon 0 + r \rightarrow r $ and $ \varepsilon\colon r\rightarrow 0 + r $ are the special congruences.

For the second equality we compute both sides. For the left side we have
\begin{equation*}
F_1 \left( \beta'(r''', c') \alpha' \circ \beta(r, c)\alpha \right) = F_1 \left( \theta(r, c' + c) \alpha\right),
\end{equation*}
where we have the following notation
\begin{equation*}
r' \xrightarrow{~~\alpha~~} r \xrightarrow{~(r,c)~} \partial(c)+r\xrightarrow{~~\beta~~} r''\xrightarrow{~~\alpha'~~} r''' \xrightarrow{~(r''',c')~} \partial(c')+r''' \xrightarrow{~~\beta'~~} \overline{r}
\end{equation*}
and $ \theta $ is a composition of special isomorphisms
\begin{equation*}
\theta\colon (\partial(c') + \partial(c)) + r \overset{\varphi}{\approx } \partial(c') +(\partial(c)+ r)  \approx \partial(c') + r'' \approx \partial(c')+r''' \xrightarrow{~~\beta'~~} \overline{r}.
\end{equation*}

Therefore we obtain
\begin{equation*}
F_1\left( \theta(r, c' + c) \alpha\right)  = \theta((c' + c) + i(r)) \varepsilon \alpha,
\end{equation*}
where $\varepsilon\colon r \rightarrow 0+r$ is the special isomorphism.

Denote $\theta((c' + c) + i(r)) \varepsilon \alpha = A$.

For the right side of \eqref{eq:51} we have
\begin{equation*}
F_1(\beta'(r''', c') \alpha') \cdot F_1(\beta(r, c) \alpha) = \beta'(c' + i(r''')) \varepsilon' \alpha' \circ \beta(c + i(r)) \varepsilon \alpha
\end{equation*}
where $\varepsilon'\colon r''' \rightarrow 0+r'''$ is the special isomorphism.

Now we denote $\varepsilon' \alpha' \beta = \tilde{\theta}$ and $(c' + i(r''')) \circ \tilde{\theta}(c + i(r)) \varepsilon \alpha = B$.

Applying the special isomorphism $fg \approx (f - i d_0(f)) + g$ (see Lemma \ref{lem:compadd}) we obtain the following special isomorphisms
\begin{alignat*}{2}
B & \approx  ((c' + i(r''')) - i d_0 (c' + i(r'''))) + \tilde{\theta}(c + i(r)) \\
& \approx  c' + \tilde{\theta}(c + i(r)) \\
& \approx  c' + (c + i(r)) \\
& \approx (c' + c) + i(r).
\end{alignat*}

Since $dom A=r'=dom \beta'B$ and $codom A=r'=codom \beta'B$, applying the coherence property we obtain the desired equality $A=\beta'B$ .
\end{proof}

\begin{prop}\label{prop:56}
$(F_0,F_1)$ is a morphism of categorical groups.
\end{prop}

\begin{proof}
The proof follows from Lemmas \ref{lem:54} and \ref{lem:55}.
\end{proof}

\begin{prop}\label{prop:57}
$PF=1_{\overline{\C}}$.
\end{prop}
\begin{proof}
Obviously, $P_0 F_0 = 1_{\mathbb{G}_0}$. For any $\beta(r,c)\alpha \in \mathbb{G}_1$, we have to show that
\begin{equation*}
P_1 F_1(\beta(r,c)\alpha) = \beta(r,c)\alpha.
\end{equation*}
Here we have the following sequence of arrows:
\begin{equation*}
r' \xrightarrow{~~\alpha~~} r \xrightarrow{~(r,c)~} \partial(c)+r \xrightarrow{~\beta~} r''.
\end{equation*}
By the definition of $ F_1 $, we have
\begin{equation*}
F_1(\beta(r,c)\alpha) = \beta(c+i(r)) \varepsilon \alpha,
\end{equation*}
where $ \varepsilon\colon r \rightarrow 0 + r $ is the special congruence in $\Str$ and $ c+i(r) $ is an arrow
\begin{equation*}
0+r\rightarrow \partial(c)+r.
\end{equation*}

By the definition of $ P_1 $, we have
\begin{alignat*}{2}
P_1 F_1(\beta(r,c)\alpha) & = P_1(\beta(c+i(r)) \varepsilon\alpha) \\
& = \theta( d_0(\beta(c+i(r))\varepsilon\alpha) , ((\beta(c+i(r))\varepsilon\alpha-i d_0(\beta(c+i(r))\varepsilon\alpha)\alpha_1)).
\end{alignat*}
where $ \varepsilon$, $\alpha_1$ and $\theta$ are the usual special congruences. We need to show that
\begin{equation*}
\theta( d_0(\beta(c+i(r))\varepsilon\alpha) , ((\beta(c+i(r))\varepsilon\alpha-i d_0(\beta(c+i(r))\varepsilon\alpha)\alpha_1)))=\beta(r,c)\alpha.
\end{equation*}

We have the following special congruence relations
\begin{alignat*}{3}
r'= d_0(\beta(c+i(r))\varepsilon\alpha) & \stackrel{\alpha}{\sim} & r \\
\beta(c+i(r))\varepsilon\alpha-i d_0(\beta(c+i(r))\varepsilon\alpha & \sim & c
\end{alignat*}

By Corollary \ref{cor:2}, we obtain the desired equality.
\end{proof}

\begin{prop}\label{prop:58}
$FP=1_{\C}$.
\end{prop}
\begin{proof}
We only need to prove that, for any arrow $ f \in \C $, $ F(P(f)) = f $.

We have
\begin{equation*}
F P(f) = F\left(\theta\left( d_0(f), \left(f - i d_0(f)\right) \alpha\right) \right) = \theta\left(\left(f - i d_0(f)\right)\alpha + i(d_0(f))\right)\varepsilon.
\end{equation*}

The picture for the arrows is the following:
\begin{equation*}
0 \xrightarrow{~~\alpha~~} d_0(f) - d_0(f) \xrightarrow{~f - i d_0(f)~} d_1(f) - d_0(f)
\end{equation*}
\begin{equation*}
d_0(f) \xrightarrow{~~\varepsilon~~} 0 + d_0(f) \xrightarrow{~\left(f - i d_0(f)\right) \alpha + i\left(d_0(f)\right)~}  \left(d_1(f) - d_0(f)\right) + d_0(f)  \stackrel{\theta}{\approx} d_1(f).
\end{equation*}

We have the special isomorphism $(f - i d_0(f))\alpha + id_0(f) \approx f$, which by the coherence property, gives the desired equality.
\end{proof}
\begin{thm}\label{theo:59}
We have an isomorphism $\mathbb{TL}\approx 1_{\Catgr}$.
\end{thm}
\begin{proof}
It follows from Propositions \ref{prop:57} and \ref{prop:58}.
\end{proof}
	
\section{The isomorphism $\mathbb{LT}\approx 1_{\cCrsM}$ and the main theorem}

Let $\partial\colon M \rightarrow N \in \cCrsM$. According to our notation we have $ \mathbb{T}\left(M \xrightarrow{~\partial~} N\right) = \overline{\C} = \left( \mathbb{G}_0, \mathbb{G}_1, \overline{d_0}, \overline{d_1}, \overline{i}, \overline{m}\right) $ and $ \mathbb{L}\left( \overline{\C}\right) = \Strb \xrightarrow{~\overline{\partial}~} N$. We will construct maps $\xymatrix{\left(M \xrightarrow{~\partial~} N\right) \ar@<0.5ex>[r]^-{\varphi} & \left( \Strb \xrightarrow{~\overline{\partial}~} N \right) \ar@<0.5ex>[l]^-{\psi} }$, where $\varphi = (\varphi_{1}, \varphi_{0})$, $\psi = (\psi_{1}, \psi_{0})$, and prove that they are c-crossed module morphisms with the properties $\psi\varphi = 1_{M \xrightarrow{~\partial~} N}$, $ \varphi \psi = 1_{\Strb \xrightarrow{~\overline{\partial}~} N}$.

For any $r \in N$ and $c \in M $ define $\varphi_{0}(r) = r$ and $\varphi_{1}(c) = \beta(0,c)$, where $\beta\colon \partial(c) + 0 \rightarrow \partial(c)$ is the special congruence.

\begin{prop}\label{prop:61}
$(\varphi_{1}, \varphi_{0})$ is a c-crossed module morphism.
\end{prop}

\begin{proof}
Since $ M $ and $ \Strb $ are connected, $\varphi_{1}$ preserves the congruence relation in $ M$. Now we will show that $\varphi_{1} $ preserves sums, i.e. for any $ c_{1}, c_{2} \in M$, $\varphi_{1}(c_{1} + c_{2}) = \varphi_{1}(c_{1}) + \varphi_{1}(c_{2})$.

We have $ \varphi_1 (c_1 + c_2) = \alpha(0, c_1+c_2) $ and $ \varphi(c_1) + \varphi(c_2) = \alpha_1(0, c_1) + \alpha_2 (0, c_2) $, where
\begin{alignat*}{2}
\alpha\colon & \partial(c_1+c_2) + 0 \rightarrow \partial(c_1) + \partial(c_2), \\
\alpha_1\colon & \partial(c_1)+ 0 \rightarrow \partial(c_1), \\
\alpha_2\colon & \partial(c_2) + 0 \rightarrow \partial(c_2)
\end{alignat*}
are the special congruences. Moreover,
\begin{equation*}
\alpha_1(0, c_1) + \alpha_2(0, c_2) = \left( \alpha_1 + \alpha_2 \right)\overline{\psi}(0 + 0, c_1 + 0\cdot c_2)\overline{\varphi},
\end{equation*}
where $ \overline{\varphi} $ and $ \overline{\psi} $ are the usual congruence relations involved in the definition of sum in $\Strb$. We have the following sequence of arrows:

\begin{equation*}
\xymatrix@C=25mm{
0 \ar[d]_-{\overline{\varphi}} & & \partial(c_1) + \partial(c_2) \\
0+0 \ar[r]_-{(0+0, c_1+0\cdot c_2)} & (\partial(c_1)+\partial(0\cdot c_2))+(0+0) \ar[r] & (\partial(c_1)+0)+(\partial(c_2)+0)  \ar[u]_-{\alpha_1+\alpha_2}
}
\end{equation*}

We have to show
\begin{equation}
(\alpha_1 + \alpha_2) \overline{\psi}(0 + 0, c_1 + 0\cdot c_2) \overline{\varphi} = \alpha(0, c_1+c_2).
\end{equation}

We have the following special congruence relations:
\begin{alignat*}{2}
0 + 0 \sim & 0, \\
c_1+0\cdot c_2 \sim & c_1 + c_2,
\end{alignat*}
from which by Corollary \ref{cor:2} follows equality (6.1).

Now we shall prove that for any $ c \in M $ and $ r \in N $
\begin{equation}\label{eq:62}
\varphi_1(r \cdot c) = r \cdot \varphi_1(c).
\end{equation}

We have $ \varphi_1(c) = \alpha(0,c) $ and $ \varphi_1(r \cdot c) = \beta(0,r \cdot c) $, where $ \alpha\colon \partial(c)+0 \xrightarrow{~\sim~} \partial(c)$ and $ \beta \colon \partial(r\cdot c)+0 \xrightarrow{~\sim~} \partial(r\cdot c) $ are special congruences. By the definition of action of $ N $ on $\Strb$, we have
\begin{equation*}
r\cdot (\alpha(0,c))=\left(\overline{i}(r)+(\alpha(0,c)- \overline{i}(r))\right)\gamma
\end{equation*}
\begin{equation*}
\overline{i}(r)=\varepsilon(r,0)
\end{equation*}
and
\begin{alignat*}{2}
-\overline{i}(r) & = - (\varepsilon(r,0)) \\
& = (-\varepsilon) \varepsilon_1 (-r, -r\cdot 0) \\
& = (-\varepsilon)\varepsilon_1(-r,0),
\end{alignat*}
where $\varepsilon\colon \partial(0)+r\longrightarrow r$, $\varepsilon_1\colon \left(-r+\left(\partial(0)-(-r)\right)\right)+(-r)\longrightarrow -r+\partial(0)$ and $\gamma\colon 0 \rightarrow r + (0-r)$ are the special isomorphisms.

Therefore we obtain
\begin{alignat*}{2}
r\cdot \left(\alpha(0,c)\right) & = \left(\varepsilon\left(r,0\right)+\left(\alpha+\left(-\varepsilon\right)\varepsilon_1\right)\theta\left(0-r,c+0\cdot\left(-r\cdot 0\right)\right)\right)\gamma \\
& = \left(\varepsilon+\left(\alpha+\left(-\varepsilon\right)\varepsilon_1\right)\theta\right)\theta'\cdot \left(r+\left(0-r\right),0+r\cdot c\right)\gamma,
\end{alignat*}
where $\theta$ and $\theta'$ are the appropriate special isomorphisms and $ \varepsilon\colon\partial(0)+r\xrightarrow{~\sim~} r$, $ \alpha\colon\bullet\xrightarrow{~\sim~} \partial(c)$ and $-\varepsilon\colon\bullet\xrightarrow{~\sim~} -r$ are the obvious special isomorphisms.

We have $\operatorname{dom} \gamma = 0 $ and $ \operatorname{codom}\left(\left(\varepsilon+\left( \alpha+\left(-\varepsilon\right)\varepsilon_1\right)\theta\right)\theta'\right) = r + \left( \partial(c)-r\right) $.

We have the following diagram:

\begin{equation*}
\xymatrix@C=10mm@R=20mm{
r+\left(0-r\right) \ar[rr]^-{\left(r+\left(0-r\right),0+r\cdot c\right)} & & \bullet \ar[d]^-{\left(\varepsilon+\left( \alpha+\left(-\varepsilon\right)\varepsilon_1\right)\theta\right)\theta'} \\
0 \ar[u]^-{\gamma} \ar[r]_-{\left(0,r\cdot c\right)} & \partial(r\cdot c)+0 \ar[r]_-{\beta} & \partial(r\cdot c) = r + \left(\partial(c)-r\right)
}
\end{equation*}
where $r+\left(0-r\right)\approx 0$ and $0+r\cdot c \approx r\cdot c$ are special isomorphisms. By Corollary \ref{cor:2} we obtain equality \eqref{eq:62}.
\end{proof}

Now we will define $\left(\psi_1,\psi_0\right)\colon\left(\Strb\rightarrow N\right)\rightarrow\left(M\xrightarrow{~\partial~}N\right)$. $\psi_0$ is defined as the identity morphism. For any element $\beta(r,c)\alpha\in\Strb$, where $r\in N$, $c\in M$ and $\alpha\colon 0\rightarrow r$ and $\beta\colon\partial(c)+r\rightarrow t$ are special isomorphisms we have the composition
\begin{equation*}
\partial(c)\xrightarrow{~~\varepsilon~~} \partial(c)+0 \xrightarrow{~1_{\partial(c)}+\alpha~} \partial(c)+r \xrightarrow{~~\beta~~} t.
\end{equation*}

Since the c-crossed module $\partial\colon M\rightarrow N$ is special, there exist a unique $c_1\in M$ with $\partial(c_1)=t$ and a weak special congruence $c_1\sim c$. We define $\psi_1\left(\beta\left(r,c\right)\alpha\right)=c_1$.

\begin{prop}\label{prop:62}
$ (\psi_1, \psi_0) $ is a c-crossed module morphism.
\end{prop}
\begin{proof}
$\psi_1$ preserves congruence relations, since both $\Strb$ and $ M $ are connected. Now we will show that $ \psi_1 $ is a morphism of c-groups. Let $\beta(r,c)\alpha$ and $ \beta'(r', c') \alpha' \in \Strb $. We have to prove that
\begin{equation}\label{eq:63}
\psi_1\left(\beta(r, c) \alpha + \beta'(r', c') \alpha'\right) = \psi_1\left(\beta(r, c)\alpha\right) + \psi_1\left(\beta'(r', c') \alpha'\right).
\end{equation}

We have $ \psi_1(\beta(r, c) \alpha) = c_1 $, $ \psi_1(\beta'(r', c') \alpha') = c_2 $ with $ \partial(c_1) = t $, $ \partial(c_2) = t' $, $ c_1 \sim c$, $c_2\sim c'$. Thus the right side of \eqref{eq:63} is equal to $ c_1 + c_2 $. We have
\begin{equation*}
\beta(r, c) \alpha + \beta'(r', c') \alpha'= \left(\beta+\beta'\right)\theta\left(r+r', c+r\cdot c'\right)\left(\alpha+\alpha'\right),
\end{equation*}
where the congruence $ \theta $ is defined by the addition rule. Moreover
\begin{equation*}
\partial\left(c+r\cdot c'\right)=\partial(c)+\left(r+\left(\partial(c')-r\right)\right)\sim t+t',
\end{equation*}
since $ r \sim 0 $,
\begin{equation*}
\psi_1\left(\beta(r, c) \alpha + \beta'(r',c') \alpha'\right) = \psi_1\left(\left(\beta+\beta'\right)\theta\left(r+r',c+r\cdot c'\right)\left(\alpha+\alpha'\right)\right).
\end{equation*}

We have $ \partial(c_1+c_2) = t+t' $, which proves \eqref{eq:63}. Now we have to prove that $ (\psi_1, \psi_0) $ is compatible with action, i.e.
\begin{equation}
\psi_1\left(r_0\cdot\left(\beta\left(r,c\right)\alpha\right)\right) = r_0\cdot\psi_1\left(\beta\left(r,c\right)\alpha\right),
\end{equation}
for any $r_0\in N$.

We have
\begin{equation*}
\psi_1\left(r_0\cdot\left(\beta\left(r,c\right)\alpha\right)\right) = \psi_1\left(\varepsilon(r_0, 0) + \left(\beta\left(r,c\right)\alpha -\varepsilon(r_0, 0) \right)\overline{\gamma}\right),
\end{equation*}
where $\varepsilon$ and $\overline{\gamma}$ are the special isomorphisms. Easy computations show that $ \psi_1\left(r_0\cdot\left(\beta(r, c) \alpha\right)\right)= \overline{c_1}$, where $ \partial(r_0\cdot c) \sim t $, $ \partial(\overline{c_1}) = t $ and $ \overline{c_1} \sim  r_0 \cdot c $. From which we obtain $ \overline{c_1} = r_0\cdot c_1$, where $ c_1 = \psi_1\left(\beta(r, c) \alpha\right) $.
\end{proof}

\begin{prop}\label{prop:63}
For a c-crossed module $ \partial\colon M \rightarrow N $, we have $ \psi\varphi = 1_M $.
\end{prop}
\begin{proof}
For any $ c \in M $ with $ \partial(c) = t $ by the definition of $ \psi $ and $ \varphi $ we have
\begin{equation*}
\psi\varphi\left(c\right)=\psi\left(\beta\left(0,c\right)\right)=c_1,
\end{equation*}
where $ c_1 \in M $ is a unique element with weak special congruence $ c_1\sim c $ and $ \partial(c_1) = t $, $ \beta \colon \partial(c)+0\rightarrow t $ is a special congruence and we have the sequence
\begin{equation*}
\partial(c)\sim\partial(c)+0 \xrightarrow{~~\beta~~} t
\end{equation*}
Since $c_1$ is unique element with such property, we have $c_1=c$.
\end{proof}

\begin{prop}\label{prop:64}
For the c-crossed module $\Strb\xrightarrow{~\overline{\partial}~} N$, we have
\begin{equation*}
\varphi\psi = 1_{\Strb\xrightarrow{~\overline{\partial}~} N}.
\end{equation*}
\end{prop}
\begin{proof}
For any $ \beta(r, c) \alpha \in \Strb $ with
\begin{equation*}
0 \xrightarrow{~~\alpha~~} r \xrightarrow{~(r,c)~} \partial(c)+r \xrightarrow{~~\beta~~} t,
\end{equation*}
we have $ \varphi \psi\left(\beta(r, c) \alpha\right) = \varphi(c_1) $, where $ c_1 $ is a unique element weak special congruent to $ c $, i.e., $ c_1 \sim c $, $\partial(c_1)=t$. We obtain $\varphi(c_1) = \beta_1(0, c_1)$ where $\beta_1\colon \partial(c_1) + 0 \rightarrow t$ is the special congruence. We have to prove that

\begin{equation}\label{eq:65}
\beta_1(0, c_1) = \beta(r, c) \alpha.
\end{equation}

We have the following diagram:
\begin{equation*}
\xymatrix@C=15mm@R=15mm{
r \ar[rr]^-{(r,c)} & & \partial(c)+r \ar[d]^-{\beta} \\
0 \ar[u]^-{\alpha} \ar[r]_-{(0,c_1)} & \partial(c_1)+0 \ar[r]_-{\beta_1} & t
}
\end{equation*}

Since $ \alpha\colon 0 \rightarrow r $ is a special congruence and $ c \sim c_1 $  is a weak special congruence, by Corollary \ref{cor:2} we obtain equality \eqref{eq:65}.%, which proves the proposition.
\end{proof}
\begin{thm}\label{theo:65}
We have an isomorphism of functors
\begin{equation*}
\mathbb{LT}\cong 1_{\cCrsM}.
\end{equation*}
\end{thm}

\begin{proof}
It follows from Propositions \ref{prop:61}, \ref{prop:62}, \ref{prop:63}, \ref{prop:64}.
\end{proof}

\begin{thm}\label{theo:3}
We have an equivalence of categories $\Catgr \simeq \cCrsM$.
\end{thm}

\begin{proof}
It follows from Theorems \ref{theo:59} and \ref{theo:65}.
\end{proof}

\end{document}